\documentclass[times]{elsarticle}

\usepackage{amsmath}
\usepackage{amsthm}
\usepackage{amssymb}
\usepackage{bm}
\journal{arXiv.org}

\newtheorem{thm}{Theorem}
\newproof{pf}{Proof}

\begin{document}

\begin{frontmatter}

\title{A splitting scheme to solve an equation for fractional powers of elliptic operators\tnoteref{label1}}
\tnotetext[label1]{This work was supported by the Russian Foundation for Basic Research (projects 14-01-00785, 15-01-00026).}

\author{Petr N. Vabishchevich\corref{cor1}\fnref{lab1,lab2}}
\ead{vabishchevich@gmail.com}
\cortext[cor1]{Correspondibg author.}

\address[lab1]{Nuclear Safety Institute, Russian Academy of Sciences,
              52, B. Tulskaya, 115191 Moscow, Russia}

\address[lab2]{North-Eastern Federal University,
	      58, Belinskogo, 677000 Yakutsk, Russia}

\begin{abstract}

An equation containing a fractional power of an elliptic operator of second order
is studied for Dirichlet boundary conditions. Finite difference approximations in space are employed.
The proposed numerical algorithm is based on solving an auxiliary 
Cauchy problem for a pseudo-parabolic equation. Unconditionally stable 
vector additive schemes (splitting schemes) are constructed. Numerical results 
for a model problem in a rectangle calculated using the splitting with respect to 
spatial variables are presented.

\end{abstract}


\end{frontmatter}

\section{Introduction}

Non-local applied mathematical models involving fractional derivatives in time and space 
are actively discussed \cite{baleanu2012fractional,eringen2002nonlocal,kilbas2006theory} at the present time. 
A lot of models in applied physics, biology, hydrology and finance
include both sub-diffusion (fractional in time) and supper-diffusion (fractional in space) operators. 
Supper-diffusion problems are treated as evolutionary problems with a fractional power of an elliptic operator.

For solving problems with fractional powers of elliptic operators, we can apply
finite volume or finite element methods oriented to using arbitrary
domains and irregular computational grids \cite{KnabnerAngermann2003,QuarteroniValli1994}.
A computational realization is associated with the implementation of the matrix function-vector multiplication.
For such problems,  various approaches \cite{higham2008functions} were developed.
The application of Krylov subspace methods with the Lanczos approximation
for solving systems of linear equations associated with
the fractional elliptic equations is discussed in \cite{ilic2009numerical}.
A comparative analysis of the contour integral method, the extended Krylov subspace method, 
and the preassigned poles and interpolation nodes method for solving
space-fractional reaction-diffusion equations is presented in \cite{burrage2012efficient}.
The simplest variant is based on the explicit construction of the solution using the known
eigenvalues and eigenfunctions of the elliptic operator with diagonalization of the corresponding matrix
\cite{bueno2012fourier,ilic2005numerical,ilic2006numerical}. 
Unfortunately, all these approaches demonstrates high computational costs for multidimensional problems.

We have proposed \cite{vabishchevich2014numerical} a numerical algorithm to solve
an equation for fractional powers of elliptic operators that is based on
a transition to a pseudo-parabolic equation.
For an auxiliary Cauchy problem, the standard two-level schemes are applied.
The computational algorithm is simple for practical use, robust, and applicable to solving
a wide class of problems. A small number of time steps is required to find a solution.
This computational algorithm for solving equations with fractional powers of operators
is promising for transient problems. 

In the study of difference schemes for time-dependent problems of
mathematical physics, the general theory of
stability (well-posedness) for operator-difference schemes
\cite{Samarskii1989,SamarskiiMatusVabischevich2002}
is in common use. At the present time,
the exact (matching necessary and sufficient) conditions for stability are obtained for
a wide class of two- and three-level difference schemes considered in
finite-dimensional Hilbert spaces. We emphasize a constructive nature of the general theory of stability
for operator-difference schemes, where stability criteria are
formulated in the form of operator inequalities, which are easy to verify.

In numerical solving initial-boundary value problems for multidimensional PDEs,
great attention is paid to the construction of additive
schemes \cite{Marchuk1990,Vabishchevich2014}. The transition to a chain of simpler problems allows us 
to construct efficient difference schemes. Here we speak of splitting with respect to the spatial variables. In some
cases, it is useful to separate subproblems of distinct nature. In this case, we say about
splitting into physical processes. Such schemes appear in
the solution of unsteady problems for systems of interconnected equations.
There are actively discussed regionally additive schemes (domain decomposition methods), which are
focused on designing computational algorithms for parallel
computers. Iterative methods for solving steady-state problems are often
treated as pseudo-time evolution methods for solving time-dependent problems.
Many iterative methods may be associated with the use of certain additive schemes. 

In this paper, we consider issues of constructing 
unconditionally stable schemes for approximate solving problems with fractional 
powers of elliptic operators on the basis of a pseudo-time evolutionary problem. 
The success is achieved through the use of  vector additive difference schemes of
multicomponent splitting \cite{Abrashin1990,Vabischevich1996b}. The original problem is reformulated as
a vector problem. In this case, instead of a single approximate solution, we search a vector
of approximate solutions. The corresponding additive schemes are schemes
of full approximation, i.e., at each time level, we search the approximate solution of the problem. 

The paper is organized as follows.
The formulation of a steady-state problem for a space-fractional elliptic equation is given in Section \ref{sec:2}.
Finite difference approximations in space and the standard two-level schemes 
to solve an auxiliary Cauchy problem for a pseudo-parabolic equation are discussed in Section \ref{sec:3}.
In Section \ref{sec:4}, we construct a special additive difference scheme for time-stepping
and investigate its stability. The results of numerical experiments are described in Section \ref{sec:5}.

\section{Problem formulation}\label{sec:2}

In a bounded polygonal domain $\Omega \subset R^2$ with the Lipschitz continuous boundary $\partial\Omega$,
we search the solution for a  problem with a fractional power of an elliptic operator.
Introduce the elliptic operator as
\begin{equation}\label{2.1}
  \mathcal{A}  u = - {\rm div}  k({\bm x}) {\rm grad} \, u + c({\bm x}) u
\end{equation} 
with coefficients $0 < k_1 \leq k({\bm x}) \leq k_2$, $c({\bm x}) \geq 0$.
The operator $\mathcal{A}$ is defined on the set of functions $u({\bm x})$ that satisfy
on the boundary $\partial\Omega$ the following conditions:
\begin{equation}\label{2.2}
  u ({\bm x}) = 0,
  \quad {\bm x} \in \partial \Omega .
\end{equation} 

In the Hilbert space $\mathcal{H} = L_2(\Omega)$, we define 
the scalar product and norm in the standard way:
\[
  (u,w) = \int_{\Omega} u({\bm x}) w({\bm x}) d{\bm x},
  \quad \|u\| = (u,u)^{1/2} .
\] 
In the spectral problem
\[
 \mathcal{A}  \varphi_k = \lambda_k \varphi_k, 
 \quad \bm x \in \Omega , 
\] 
\[
  \varphi_k = 0,
  \quad {\bm x} \in \partial \Omega , 
\] 
we have 
\[
 \lambda_1 \leq \lambda_2 \leq ... ,
\] 
and the eigenfunctions  $ \varphi_k, \ \|\varphi_k\| = 1, \ k = 1,2, ...  $ form a basis in $L_2(\Omega)$. Therefore, 
\[
 u = \sum_{k=1}^{\infty} (u,\varphi_k) \varphi_k .
\] 
Let the operator $\mathcal{A}$ be defined in the following domain:
\[
 D(\mathcal{A} ) = \Big \{ u \ | \ u(x) \in L_2(\Omega), \ \sum_{k=0}^{\infty} | (u,\varphi_k) |^2 \lambda_k < \infty \Big \} .
\] 
Under these conditions  $\mathcal{A} : L_2(\Omega) \rightarrow L_2(\Omega)$ and
the operator $\mathcal{A}$ is self-adjoint and positive definite: 
\begin{equation}\label{2.3}
  \mathcal{A}  = \mathcal{A} ^* \geq \lambda_1 \mathcal{I} ,
  \quad \lambda_1 > 0 ,    
\end{equation} 
where $\mathcal{I}$ is the identity operator in $\mathcal{H}$.
In applications, the value of $\lambda_1$ is unknown (the spectral problem must be solved).
Therefore, we suppose that $\delta \leq \lambda_1$ in (\ref{2.3}).
Let us assume for the fractional power of the  operator $\mathcal{A}$:
\[
 \mathcal{A} ^\alpha u =  \sum_{k=0}^{\infty} (u,\varphi_k) \lambda_k^\alpha  \varphi_k .
\] 
More general and mathematically complete definition of fractional powers of elliptic operators 
is given in \cite{yagi2009abstract}. 
The solution $u(\bm x)$ satisfies the equation
\begin{equation}\label{2.4}
  \mathcal{A}^\alpha u = f 
\end{equation} 
under the restriction $0 < \alpha < 1$.

\section{Numerical algorithm}\label{sec:3}

We consider the simplest case, where the computational domain $\Omega$ is a rectangle:
\[
 \Omega = \{ \bm x  \ | \ \bm x = (x_1,x_2), \ 0 < x_k < l_k, \ k = 1,2 \} .
\]
To solve approximately the problem  (\ref{2.4}), we introduce in the domain $\Omega$ a uniform grid
\[
\overline{\omega}  = \{ \bm{x} \ | \ \bm{x} =\left(x_1, x_2\right), \quad x_k =
i_k h_k, \quad i_k = 0,1,...,N_k,
\quad N_k h_k = l_k, \ k = 1,2 \} ,
\]
where $\overline{\omega} = \omega \cup \partial \omega$ and
$\omega$ is the set of interior nodes, whereas $\partial \omega$ is the set of boundary nodes of the grid.
For grid functions $y(\bm x)$ such that $y(\bm x) = 0, \ \bm x \notin \omega$, we define the Hilbert space
$H=L_2\left(\omega\right)$, where the scalar product and the norm are specified as follows:
\[
\left(y, w\right) =  \sum_{\bm x \in  \omega} y\left(\bm{x}\right)
w\left(\bm{x}\right) h_1 h_2,  \quad 
\| y \| =  \left(y, y\right)^{1/2}.
\]

For the discrete operator $A$, we  use the additive representation
\[
A = \sum_{k=1}^{2} A_k, 
\quad \bm{x} \in \omega,
\]
where $A_k, \ k=1,2$ are associated with the corresponding differential operator with the second derivative
in one direction.

For all grid points except adjacent to the boundary, the grid operator $A_1$  can be written as
\[
  \begin{split}
  A_1 y = & -
  \frac{1}{h_1^2} k(x_1+0.5h_1,h_2) (y(x_1+h_1,h_2) - y(\bm{x})) \\ 
  & + \frac{1}{h_1^2} k(x_1-0.5h_1,h_2) (y(\bm{x}) - y(x_1-h_1,h_2)) + c(\bm x) y(\bm x), \\
  & \qquad \bm{x} \in \omega, 
  \quad x_1 \neq 0.5h_1,
  \quad x_1 \neq l_1-0.5h_1. 
 \end{split} 
\]
At the points that are adjacent to the boundary, the approximation is constructed taking into account
the boundary condition (\ref{2.2}):
\[
  \begin{split}
  A_1 y = & -
  \frac{1}{h_1^2} k(x_1+0.5h_1,h_2) (y(x_1+h_1,h_2) - y(\bm{x}))  \\ 
  & + \frac{1}{h_1^2} k(x_1-0.5h_1,h_2) y(\bm{x}) + c(\bm x) y(\bm x), \\
  & \qquad \bm{x} \in \omega, 
  \quad x_1 = 0.5h_1,
 \end{split} 
\] 
\[
  \begin{split}
  A_1 y = & 
  \frac{1}{h_1^2} k(x_1+0.5h_1,h_2)  y(\bm{x})  \\ 
  + & \frac{1}{h_1^2} k(x_1-0.5h_1,h_2) (y(\bm{x}) - y(x_1-h_1,h_2))  + c(\bm x) y(\bm x), \\
  & \qquad \bm{x} \in \omega, 
  \quad x_1 = l_1-0.5h_1. 
 \end{split} 
\]
Similarly, we construct the grid operator $A_2$.
For the above grid operators (see, e.g., \cite{Samarskii1989,SamarskiiNikolaev1978}), we have
\[
  A_k = A_k^* \geq \delta_k I, 
  \quad  \delta_k > 0 , 
  \quad  k=1,2 . 
\]
where $I$ is the grid identity operator. Because of this, the discrete operator $A$
is self-adjoint and positive definite in $H$:
\begin{equation} \label{3.1}
 A = A^* \geq \delta I,
 \quad \delta =  \sum_{k=1}^{2} \delta_k.
\end{equation}
For problems with sufficiently smooth coefficients and the right-hand side, 
it approximates the differential operator with the truncation error 
$\mathcal{O} \left(|h|^2\right)$, $|h|^2 = h_1^2+h_2^2$. 

To handle the fractional power of the grid operator $A$, let us consider the eigenvalue problem
\[
 A \varphi_m = \lambda_m^h\varphi_m . 
\] 
We have
\[
 \delta = \lambda_1^h \leq \lambda_2^h   \leq ... \leq \lambda_M^h,
 \quad M = (N_1-1)(N_2-1) , 
\] 
where eigenfunctions $\varphi_m, \ \|\varphi_m\| = 1, \ m = 1,2, ..., M.$ form a basis in $H$. Therefore
\[
 y = \sum_{m= 1}^{M}(y, \varphi_m) \varphi_m . 
\]
For the fractional power of the operator $A$, we have
\[
 A^\alpha y = \sum_{m= 1}^{M}(y, \varphi_m) (\lambda_m^h)^\alpha \varphi_m .
\] 

Using the above approximations, we arrive from (\ref{2.4}) at the discrete problem
\begin{equation}\label{3.2}
 A^\alpha w = f .
\end{equation} 

An approximate solution is sought as a solution of an auxiliary pseudo-time evolutionary problem 
\cite{vabishchevich2014numerical}.
Assume that
\[
 y(t) = (\theta \delta)^{\alpha} (t (A - \theta \delta I) + \theta  \delta I)^{-\alpha} y(0) ,
\]
with $0 < \theta < 1$. Therefore
\[
 y(1) =  (\theta \delta)^{\alpha} A^{-\alpha} y(0)
\]  
and then $w = y(1)$.
The function $y(t)$ satisfies the evolutionary equation
\begin{equation}\label{3.3}
  (t D + \theta \delta I) \frac{d y}{d t} + \alpha D y = 0 ,
  \quad 0 < t \leq 1 ,
\end{equation}  
where
\[
 D = A - \theta \delta I .
\] 
By (\ref{3.1}), we get
\begin{equation}\label{3.4}
 D = D^* \geq (1-\theta) \delta I > 0 .
\end{equation} 
We supplement (\ref{3.3}) with the initial condition
\begin{equation}\label{3.5}
 y(0) = (\theta \delta)^{-\alpha} f .  
\end{equation} 
The solution of equation (\ref{2.4}) can be defined as the solution of the Cauchy problem 
(\ref{3.3})--(\ref{3.5}) at the final pseudo-time moment $t=1$.
In \cite{vabishchevich2014numerical}, the case with $\theta =1$ was considered.

For the solution of the problem (\ref{3.3}), (\ref{3.5}), it is possible to obtain various a priori estimates.
Elementary estimates have the form
\begin{equation}\label{3.6}
  \|y(t)\|_G \leq \|y(0)\|_G , 
\end{equation} 
where, e.g., $G = I, D$.
To get (\ref{3.6}) for  $G = D$, multiply scalarly equation (\ref{3.3}) by $dy/ dt$.
If $G = I$, then equation  (\ref{3.3}) is multiplied by $\alpha y + t dy/ dt$.

To solve numerically the problem (\ref{3.3}), (\ref{3.5}),
we use the simplest implicit two-level scheme.
Let $\tau$ be a step of a uniform grid in time such that $y^n = y(t^n), \ t^n = n \tau$, $n = 0,1, ..., N, \ N\tau = 1$.
Let us approximate equation (\ref{3.3}) by the implicit two-level scheme
\begin{equation}\label{3.7}
 (t^{\sigma(n)} D + \theta \delta I) \frac{ y^{n+1} - y^{n}}{\tau }
 + \alpha D y^{\sigma(n)} = 0,  \quad n = 0,1, ..., N-1,
\end{equation}
\begin{equation}\label{3.8}
 y^0 = (\theta \delta)^{-\alpha} f .
\end{equation} 
We use the notation
\[
  t^{\sigma(n)} = \sigma t^{n+1} + (1-\sigma) t^{n},
  \quad y^{\sigma(n)} = \sigma y^{n+1} + (1-\sigma) y^{n}.
\]
For $\sigma =0.5$, the difference scheme (\ref{3.7}), (\ref{3.8}) approximates the problem 
(\ref{3.3}), (\ref{3.5})
with the second order by $\tau$, whereas for other values of $\sigma$, we have only the first order.

We have
\[
 y^{\sigma(n)} = \left( \sigma - \frac{1}{2} \right ) \tau \frac{y^{n+1} - y^{n}}{\tau }   +
 \frac{1}{2} (y^{n+1} + y^{n}) .
\] 
\begin{thm}\label{t-1}
For $\sigma \geq 0.5$ and $0 < \theta < 1$, the difference scheme (\ref{3.7}), (\ref{3.8}) 
is unconditionally stable with respect to the initial data.
The approximate solution satisfies the estimate 
\begin{equation}\label{3.9}
  \|y^{n+1}\|_G \leq \|y^0\|_G , 
  \quad n = 0,1, ..., N-1,
\end{equation} 
with $G = I, D$.
\end{thm} 
 
\begin{pf}
Rewrite equation (\ref{3.7}) as
\[
 \left(t^{\sigma(n)} D + \theta \delta I + \left( \sigma - \frac{1}{2} \right ) \tau D \right )  \frac{y^{n+1} - y^{n}}{\tau } 
 + \frac{\alpha }{2} D (y^{n+1} + y^{n}) = 0 .
\] 
Multiplying scalarly this equation by $y^{n+1} - y^{n}$, for $\sigma \geq 0.5$, we get
\[
  \|y^{n+1}\|_D \leq \|y^n\|_D , 
  \quad n = 0,1, ..., N-1 .  
\] 
This inequality ensures the estimate (\ref{3.9}) for  $G = D$.

In a similar way, we consider the case with $G = I$.
Rewrite equation (\ref{3.7}) in the following form:
\[
 \theta \delta \frac{ y^{n+1} - y^{n}}{\tau } + D \left( \alpha y^{\sigma(n)} + t^{\sigma(n)}\frac{ y^{n+1} - y^{n}}{\tau } \right ) = 0 . 
\] 
Multiplying scalarly it by
\[
 \alpha y^{\sigma(n)} + t^{\sigma(n)}\frac{ y^{n+1} - y^{n}}{\tau },
\] 
in view of (\ref{3.4}), we arrive at
\[
 \left ( \frac{ y^{n+1} - y^{n}}{\tau }, y^{\sigma(n)} \right ) \leq 0 . 
\] 
If $\sigma \geq 0.5$, then
\[
  \|y^{n+1}\| \leq \|y^n\| , 
  \quad n = 0,1, ..., N-1 .  
\] 
Thus, we obtain (\ref{3.9}) for $G = I$.
\end{pf}

\section{Splitting schemes}\label{sec:4}

The numerical implementation of the difference scheme (\ref{3.7}), (\ref{3.8}) involves
the solution of standard elliptic boundary value problems
\[
 (t^{\sigma(n)} + \alpha \sigma \tau) D y^{n+1} + \theta \delta y^{n+1} = \varphi^{n} 
\] 
with the given $\varphi^n$ for $n = 0,1, ..., N-1$. 

The inversion of the operator $\theta \delta I + (t^{\sigma(n)} + \alpha \sigma \tau) (A - \theta \delta I) $ 
may be enough difficult.
Thus, it seems natural to construct difference schemes for unsteady
problems such that they will be unconditionally stable, but at the same time, their implementation would be
considerably simpler. The most interesting
results have been obtained taking into account a special structure of the problem operator $A$.

We define a class of the problems (\ref{3.1}), (\ref{3.2}), where the operator $A$
has the following $p$-component additive representation:
\begin{equation}\label{4.1}
  A = \sum_{i = 1}^{p} A_i.
\end{equation}
Assume that the operators $A_{i}, \ i=1,2,\ldots,p$ are 
simpler than $A$. We organize computations in such a way that the transition to
a new time level in solving the problem (\ref{3.3}), (\ref{3.5})
is not more complicated than the solution of the $p$ problems
for the individual operator terms:
$A_{i}, \ i=1,2,\ldots,p$. 
Such difference schemes are called additive difference schemes \cite{Vabishchevich2014}. 

A decomposition into operator terms in the additive representation (\ref{4.1}) may have a different 
nature. In particular, additive difference schemes are used to solve numerically
multidimensional transient problems of mathematical physics, where
one-dimensional problems are the most simple ones. On the basis of splitting
with respect to the spatial variables, we construct locally one-dimensional schemes.
In using computational algorithms of domain decomposition, which focus on modern parallel computers,
the original problem is divided into several subproblems, and each of them is solved in its own subdomain
on its individual processor. 

To solve numerically the evolutionary problem (\ref{3.3})--(\ref{3.5}), it seems reasonable
to apply the following additive representation of $D$:
\begin{equation}\label{4.2}
 D = \sum_{i = 1}^{p} D_i, 
 \quad D_i = D^*_i \geq  \mu_i I,
 \quad  \mu_i > 0 ,
 \quad  i=1,2,\ldots,p ,
\end{equation} 
with pairwise non-permutable operators $D_i, \ i=1,2,\ldots,p$:
\[
 D_i D_j \neq D_j D_i, 
 \quad i \neq j,
 \quad  i,j =1,2,\ldots,p .
\] 
For the splitting (\ref{4.1}), we put
\[
 D_i = A_i - \chi_i I,
 \quad  i=1,2,\ldots,p ,
 \quad \sum_{i = 1}^{p} \chi_i = \theta \delta 
\] 
with an appropriate choice of the constants $\chi_i, \ i=1,2,\ldots,p$.
In particular, if
\[
 A_i \geq \delta_i I, 
 \quad \delta_i > 0, 
\] 
then we may take
\[
 \chi_i = \theta \delta_i,
 \quad  i =1,2,\ldots,p .
\] 

The problem (\ref{3.3}), (\ref{3.5}), (\ref{4.2}) does not allow the use of the standard splitting schemes because
here we have a splitting of the operator at the time derivative. 
In the work \cite{vabishchevich2012new}, we proposed and investigated vector additive-operator schemes 
with a splitting of the operator at the time derivative into the sum of positive definite self-adjoint operators. 
Here we consider a class of problems with a splitting of both the leading operator of the problem and the operator 
at the time derivative. Such problems are typical in studying  boundary value problems for pseudo-parabolic equations. 
Some vector splitting schemes for pseudo-parabolic equations with constant operators 
are considered in \cite{vabishchevich2013splitting}.

For the Cauchy problem (\ref{3.3}), (\ref{3.5}), (\ref{4.2}), we specify the vector
${\bm y}  =  \{y_1, y_2, ..., y_p \}$.
Each individual component is determined as the solution of the similar problems
\begin{equation}\label{4.3}
   \theta \delta \frac{d y_{i}}{d t} + 
   t \sum_{j  = 1} ^{p} D_j  \frac {d y_{j }} {d t}
   + \alpha \sum_{j  = 1} ^{p} D_j  y_{j }  = 0,
   \quad 0 < t \leq 1 ,
\end{equation}
\begin{equation}\label{4.4}
   y_{i}(0) =  (\theta \delta)^{-\alpha} f ,
  \quad i =1,2,...,p .
\end{equation}

Here is an elementary coordinate-wise estimate for the solution stability. 
Subtracting one equation from another, we get
\[
  \frac{d y_{i}}{d t} - \frac{d y_{i-1}}{d t}  = 0,
  \quad i =2,3,...,p .
\]
In view of the initial conditions (\ref{4.4}), we have
\[
  y_{i} = y_{i-1} ,
  \quad i =2,3,...,p .
\]
For the individual component $y_{i}$, we obtain the same equation as for $y$:
\[
   \theta \delta \frac{d y_{i}}{d t} + 
   t \sum_{j  = 1} ^{p} D_j  \frac {d y_{i}} {d t}
   + \alpha \sum_{j  = 1} ^{p} D_j  u_{i }  = 0,
   \quad 0 < t \leq 1 ,
  \quad i =1,2,...,p .
\]
Thus (see (\ref{3.6})), a priori estimates
\begin{equation}\label{4.5}
  \|y_{i}(t)\|_G \le (\theta \delta)^{-\alpha} \|f\|_G,
  \quad i =1,2,...,p 
\end{equation}
are satisfied. Because of this, we have
\[
  y_{i}(t) = y(t) ,
  \quad 0 < t \leq 1 ,
  \quad i =1,2,...,p .
\]
and so, we can treat any component of the vector ${\bm y}(t)$
as the solution of the original problem (\ref{3.3}), (\ref{3.5}).

For the vector evolutionary problem (\ref{4.3}), (\ref{4.4}), it is easy to obtain
(see, e.g., \cite{Vabishchevich2014}) a priori estimates for the vector ${\bm y}$.
Let us introduce the Hilbert space $ {\bm H} = H^p$ with the scalar product
\[
  ({\bm y}, {\bm v}) =
  \sum_{i=1} ^{p} (y_i, v_i) .
\]

Rewrite equation (\ref{4.3}) as
\[
   \theta \delta D_i \frac{d y_{i}}{d t}+ 
   t \sum_{j  = 1} ^{p} D_i D_j  \frac {d y_{j }} {d t}
   + \alpha \sum_{j  = 1} ^{p} D_i D_j  u_{j } = 0,
  \quad 0 < t \leq 1 ,
  \quad \alpha =1,2,...,p .
\]
This allows us to write the system of equations in the vector form
\begin{equation}\label{4.6}
  {\bm B} \frac{d {\bm y}}{d t} +
  {\bm A} {\bm y} = 0.
\end{equation}
Operator matrices ${\bm B}$ and ${\bm A}$ seem like this
\begin{equation}\label{4.7}
\begin{split}
  {\bm B} & = \{ B_{i j } \},
  \quad B_{i j } = \delta_{i j }  \theta \delta  D_i + t D_i D_j, \\
  {\bm A} & = \{ A_{i j } \},
  \quad A_{i j } = \alpha D_i D_j,
  \quad i, j  =1,2,...,p ,
\end{split}
\end{equation}
where $\delta_{i j }$ is the Kronecker delta.
The equation (\ref{4.6}) is supplemented with the initial condition
\begin{equation}\label{4.8}
  {\bm y}(0) = {\bm y}^0,
\end{equation}
where ${\bm y}^0 = \{y_1(0), y_2(0), ..., y_p(0)\}$.

The principal advantage of using the formulation (\ref{4.6}) results from the fact that in ${\bm H}$,
we have
\[
  {\bm B} = {\bm B}^* > 0,
  \quad {\bm A} = {\bm A}^* \geq  0.
\]
Let us consider an a priori estimate for the solution of the vector problem (\ref{4.6})--(\ref{4.8}),
which, on the one hand, is more complicated than (\ref{4.5}), and on the other hand, 
it will serve us as a guide for constructing operator-difference schemes.

For ${\bm B}$, we introduce the representation
\[
 {\bm B} = {\bm C} + \frac{t}{\alpha} {\bm A},
\] 
where
\[
  {\bm C} = \{ C_{i j } \},
  \quad C_{i j } = \delta_{i j }  \theta \delta  D_i,
  \quad i, j  =1,2,...,p .  
\] 
Rewrite equation (\ref{4.6}) in the form
\begin{equation}\label{4.9}
  {\bm C} \frac{d {\bm y}}{d t} + 
  {\bm A} \left (  \frac{t}{\alpha}\frac{d {\bm y}}{d t} + {\bm y} \right ) = 0 . 
\end{equation} 
By (\ref{4.7}), we get
\[
  ({\bm A} {\bm y}, {\bm y} ) 
  = \alpha \left (\sum_{j=1}^{p} A_j y_{j},
  \sum_{j=1}^{p} A_j y_{j} \right ) .
\]
and so we have ${\bm A} \geq 0$. Multiply (\ref{4.9}) scalarly in ${\bm H}$ by
\[
  \frac{t}{\alpha}\frac{d {\bm y}}{d t} + {\bm y} .
\] 
This gives
\begin{equation}\label{4.10}
  \frac{t}{\alpha} \left ({\bm C} \frac{d {\bm y}}{d t},  \frac{d {\bm y}}{d t} \right ) 
  + \frac{1}{2} \frac{d}{d t} ({\bm C} {\bm y}, {\bm y}) = 0 .
\end{equation}  
In view of ${\bm C} > 0$, from (\ref{4.10}), we get that the estimate
\begin{equation}\label{4.11}
  \| {\bm y} \|^2_{{\bm C}} \leq 
  \| {\bm y}^0 \|^2_{{\bm C}}
\end{equation}
holds. Taking into account the representation for ${\bm C}$, we have
\[
  \| {\bm y} \|^2_{{\bm C}} = 
  \theta \delta \sum_{i =1}^{p}\left (  D_i y_i, y_i \right ) .
\]
Thus, the estimate (\ref{4.11}) along with (\ref{4.5}) may be treated 
as a vector analogue of the estimate (\ref{3.6}). In view of (\ref{4.2}), 
the estimate (\ref{4.11}) ensures stability of each individual 
components of the vector ${\bm y}(t)$.

To construct splitting schemes for solving the problem (\ref{3.3}), (\ref{3.5}),
we apply common schemes with weights and the vector problem (\ref{4.6}), (\ref{4.8}).
In fact, we construct splitting schemes for the system of unsteady
equations (\ref{4.3}), (\ref{4.4}), which are coupled via time derivatives. 
For such problems, the construction of unconditionally stable
splitting schemes can be carried out on the basis of the triangular splitting of the operator
matrices ${\bm B}$ and ${\bm A}$ or the separation of the diagonal part of these operator matrices
\cite{gaspar2014explicit,vabishchevich2012new,vabishchevich2014new,Vabishchevich2014}. 
In some cases, it is reasonable to combine these approaches, namely, to conduct the triangular splitting
of the operator matrix ${\bm B}$ and to separate the diagonal of the operator matrix ${\bm A}$, or vice versa.

Let us construct additive operator-difference schemes using the splitting of operator ${\bm  A}$
with separation of the diagonal part. In this case, we obtain
\begin{equation}\label{4.12}
  {\bm  A} = {\bm  A}_0 + {\bm  A}_1,
  \quad  {\bm A}_0 = \mathrm{diag} (A_{11}, A_{22}, ..., A_{pp} ) .
\end{equation} 
In additive representation (\ref{4.12}), we have
\[
  {\bm  A}_0 = 
  \begin{pmatrix}
  A_{11} & 0 & \cdots & 0 \\
  0  & A_{22} & \cdots & 0 \\
  \cdots & \cdots & \cdots & \cdots \\
  0 & 0 & \cdots &  A_{pp} \\
  \end{pmatrix} ,
  \quad
  {\bm  A}_1 = 
  \begin{pmatrix}
  0 & A_{12} & \cdots & A_{1p} \\
  A_{21} & 0 & \cdots & A_{2p} \\
  \cdots & \cdots & \cdots &  \cdots \\
  A_{p1} & A_{p2} & \cdots & \ 0 \\
  \end{pmatrix}  .
\] 
In view of
\[
\begin{split}
  \left (\sum_{j=1}^{p} A_j y_{j}, \sum_{j=1}^{p} A_j y_{j} \right ) & = 
  \left (\left ( \sum_{i =1}^{p} 
  A_j y_{j} \right )^2,  1 \right ) \\
  & \leq p \sum_{i =1}^{p}   \left ((A_j y_{j})^2, 1 \right ) 
  = p  \sum_{i =1}^{p} \left (A_j y_{j}, A_j y_{j} \right ), 
\end{split}
\] 
we get
\[
  {\bm A} \leq p {\bm A}_0 .
\]
We can consider problem (\ref{4.6}), (\ref{4.8}), (\ref{4.12}) under the additional assumption:
\begin{equation}\label{4.13}
   {\bm A}_0 \geq  \frac{1}{p} {\bm A} \geq 0 .
\end{equation} 

To solve the vector problem, we apply the following scheme with weights $\sigma_1, \sigma_2$:
\begin{equation}\label{4.14}
\begin{split}
  {\bm C} \frac{{\bm y}^{n+1} - {\bm y}^{n}}{\tau } & +
  \frac{t^{n}}{\alpha} {\bm A}_0 \left (\sigma_1 \frac{{\bm y}^{n+1} - {\bm y}^{n}}{\tau } 
  + (1-\sigma_1) \frac{{\bm y}^{n} - {\bm y}^{n-1}}{\tau } \right ) \\
  & +  \frac{t^{n}}{\alpha} {\bm A}_1 \frac{{\bm y}^{n} - {\bm y}^{n-1}}{\tau } \\
  & + {\bm A}_0 (\sigma_2 {\bm y}^{n+1} + (1-\sigma_2) {\bm y}^{n} )
  + {\bm A}_1{\bm y}^{n} = 0 , \\
  & \quad n = 1, 2, ... , N-1 .
\end{split} 
\end{equation} 
In the coordinate-wise formulation, the scheme (\ref{4.14}) corresponds to applying  the scheme
\[
\begin{split}
  \theta \delta  \frac{y_i^{n+1} - y_i^{n}}{\tau} & + 
  t^{n} D_i \left ( \sigma_1  \frac{y_i^{n+1} - y_i^{n}}{\tau} 
  + (1-\sigma_1 ) \frac{y_i^{n} - y_i^{n-1}}{\tau} \right ) 
  + t^{n} \sum_{i \neq j= 1}^{p} D_j \frac{y_j^{n} - y_j^{n-1}}{\tau } \\
  & + \alpha D_i (\sigma_2 y_i^{n+1} + (1-\sigma_2) y_i^{n} )
   + \alpha \sum_{i \neq j= 1}^{p} D_j y_j^{n} = 0, \\
  & \quad n = 0, 1, ... , N-1, 
  \quad i =1,2,...,p ,
\end{split}
\]
to the Cauchy problem (\ref{4.3}), (\ref{4.4}).
In contrast to (\ref{3.7}), the scheme (\ref{4.12}) is  three-level and 
it has the weighting parameter $\sigma_1$ in the approximation of the time derivative.

Computational implementation of (\ref{4.14}) involves the solution of grid problems
\[
  \left ( \theta \delta I + \left (\sigma_1 t^{n} + \sigma_2 \tau \alpha \right ) D_i \right ) y_i^{n+1} = 
  \chi_i^n,
  \quad i =1,2,...,p ,
\]
for the transition from the time level $t^n$ to the new level $t^{n+1}$. 

Using the notation (\ref{4.7}),  (\ref{4.12}), rewrite the operator-difference scheme (\ref{4.14})
in the canonical form of three-level schemes \cite{Samarskii1989,SamarskiiMatusVabischevich2002}:
\begin{equation}\label{4.15}
  {\bm G}^n \frac{{\bm y}^{n+1} - {\bm y}^{n-1}}{2\tau } + 
  {\bm R}^n \frac{{\bm y}^{n+1} - 2{\bm y}^{n} + {\bm y}^{n-1}}{\tau^2 }  + 
  {\bm A} {\bm y}^{n} = 0 . 
\end{equation}
In view of
\[
  \frac{{\bm y}^{n+1} - {\bm y}^{n}}{\tau } =
  \frac{{\bm y}^{n+1} - {\bm y}^{n-1}}{2 \tau }
  + \frac{{\bm y}^{n+1} - 2{\bm y}^{n+1} + {\bm y}^{n-1}}{2 \tau },
\]
\[
  \frac{{\bm y}^{n} - {\bm y}^{n-1}}{\tau } =
  \frac{{\bm y}^{n+1} - {\bm y}^{n-1}}{2 \tau }
  - \frac{{\bm y}^{n+1} - 2{\bm y}^{n+1} + {\bm y}^{n-1}}{2 \tau },
\]
for  ${\bm G}$ and  ${\bm R}$, we have
\begin{equation}\label{4.16}
  {\bm G}^n = {\bm C}  + \frac{t^{n}}{\alpha} {\bm A} + \sigma_2 \tau {\bm A}_0 ,
  \quad  {\bm R}^n = \frac{\tau }{2} {\bm C} + \sigma_1 \tau \frac{t^{n}}{\alpha} {\bm A}_0 -
  \frac{\tau }{2} \frac{t^{n}}{\alpha} {\bm A} + \sigma_2 \frac{\tau^2 }{2} {\bm A}_0 .
\end{equation}
It is essential that the operators ${\bm G}$ and ${\bm R}$ are variable, namely, they depend on time. 
When considering the three-level schemes, we have complicated norms, and this unsteadiness 
of the operators makes practically impossible to obtain global estimates for stability. 
For this reason, here we formulate a more particular result.

\begin{thm}\label{t-2}
If $2\sigma_1\geq p$ and $2\sigma_2\geq p$, then for the solution of the explicit-implicit 
scheme  (\ref{4.12}), (\ref{4.14}), the following estimate with respect to the initial data
\begin{equation}\label{4.17}
\begin{split}
  \left \| \frac{{\bm y}^{n+1} + {\bm y}^{n}}{2} 
  \right \|^2_{{\bm A}} & +
  \left \|   \frac{ {\bm y}^{n+1} - {\bm y}^{n} } {\tau} \right \|^2_{{\bm R}^n - \frac{\tau^2}{4} {\bm A}}  \\
  & \leq 
  \left \| \frac{{\bm y}^{n} + {\bm y}^{n-1}}{2} 
  \right \|^2_{{\bm A}} +
  \left \|   \frac{ {\bm y}^{n} - {\bm y}^{n-1} } {\tau} \right \|^2_{{\bm R}^n - \frac{\tau^2}{4} {\bm A}} 
\end{split}
\end{equation}
holds.
\end{thm} 
 
\begin{pf}
Rewrite (\ref{4.15}) as
\begin{equation}\label{4.18}
\begin{split}
  {\bm G}^n \frac{{\bm y}^{n+1} - {\bm y}^{n-1}}{2\tau } & + 
  \left ( {\bm R}^n - \frac{\tau^2}{4} {\bm A} \right ) \frac{{\bm y}^{n+1} - 2{\bm y}^{n} + {\bm y}^{n-1}}{\tau^2 }  \\
  & + 
  {\bm A} \frac{{\bm y}^{n+1} + 2{\bm y}^{n} + {\bm y}^{n-1}}{4}  = 0 . 
\end{split}
\end{equation}
Introduce
\[
  {\bm v}^{n} = \frac{1}{2} ({\bm y}^{n} + {\bm y}^{n-1}),
  \quad {\bm w}^{n} = \frac{{\bm y}^{n} - {\bm y}^{n-1}} {\tau} 
\]
and rewrite (\ref{4.18}) in the form
\begin{equation}\label{4.19}
  {\bm G}^n  \frac{{\bm w}^{n+1} + {\bm w}^{n}}{2 } 
  + \left ( {\bm R}^n - \frac{\tau^2}{4} {\bm A} \right ) 
  \frac{{\bm w}^{n+1} - {\bm w}^{n}}{\tau } +
  \frac{1 }{2} {\bm A}
  ( {\bm v}^{n+1} + {\bm v}^{n})  = 0.
\end{equation}
Multiplying  scalarly (\ref{4.19}) by
\[
  2 ({\bm v}^{n+1} - {\bm v}^{n}) =
  \tau ( {\bm w}^{n+1} + {\bm w}^{n} ),
\]
we obtain the equality
\begin{equation}\label{4.20}
\begin{split}
  \frac{\tau }{2 } 
  ( {\bm G}^n ({\bm w}^{n+1} + {\bm w}^{n}),
    {\bm w}^{n+1} + {\bm w}^{n}) & +
   \left ( \left ( {\bm R}^n - \frac{\tau^2}{4} {\bm A} \right )  ({\bm w}^{n+1} - {\bm w}^{n}),
    {\bm w}^{n+1} + {\bm w}^{n} \right ) \\
  & + 
  ( {\bm A} ({\bm v}^{n+1} + {\bm v}^{n}),
    {\bm v}^{n+1} - {\bm v}^{n}) = 0 .
\end{split}
\end{equation}
Taking into account the positive definiteness of the operator ${\bm G}^n$ and the self-adjointness 
of the operators ${\bm R}^n$, ${\bm A}$, from (\ref{4.20}), we get
\[
\begin{split}
 \left ( \left ( {\bm R}^n - \frac{\tau^2}{4} {\bm A} \right ) {\bm w}^{n+1}, {\bm w}^{n+1} \right ) & -
 \left ( \left ( {\bm R}^n - \frac{\tau^2}{4} {\bm A} \right ) {\bm w}^{n}, {\bm w}^{n} \right ) \\
 & +
 ( {\bm A} ({\bm v}^{n+1}, {\bm v}^{n+1}) - ( {\bm A} ({\bm v}^{n}, {\bm v}^{n}) \leq  0.
\end{split}
\]
From this inequality, we obtain (\ref{4.17}) under the condition
\begin{equation}\label{4.21}
 {\bm R}^n > \frac{\tau^2}{4} {\bm A} .
\end{equation} 
By (\ref{4.13}), (\ref{4.16}),  we get
\[
\begin{split}
 {\bm R}^n - \frac{\tau^2}{4} {\bm A} & = 
 \frac{\tau }{2} {\bm C} + 
 \frac{\tau }{2} \frac{t^{n}}{\alpha} (2 \sigma_1 {\bm A}_0 - {\bm A})  
 +  \frac{\tau^2 }{4} (2\sigma_2 {\bm A}_0 -  {\bm A}) \\
 & \geq \frac{\tau }{2} {\bm C} + \frac{\tau }{2} \frac{t^{n}}{\alpha} (2\sigma_1 - p)  {\bm A}_0 
 +  \frac{\tau^2 }{4} (2\sigma_2 - p) {\bm A}_0 > 0
\end{split}
\] 
for $2\sigma_1 - p \geq 0$ and $2\sigma_2 - p \geq 0$.
\end{pf}

\section{Numerical experiments}\label{sec:5}

For simplicity, we restrict ourselves to the case, where $\mathcal{A}$ is the Laplace operator, i.e., 
\[
 k(\bm x) = 1,
 \quad  c(\bm x) = 0 
\] 
in equation (\ref{2.1}). 
Under these assumptions, the solution of the spectral problem for the differential 
and discrete Laplace operator is well known and so, it is possible to construct the
exact solutions. In particular, for the constant $\delta$ in the inequality (\ref{3.1}), we have
\begin{equation}\label{5.1}
 \delta_k  = \frac{4}{h_k^2} \sin^2 \frac{\pi }{2 N_k} , 
  \quad  k=1,2 .  
\end{equation} 
Let the right-hand side is given as
\[
 f(\bm x) = \sin(\pi x_1)  \sin(\pi x_2) + \sin(3\pi x_1)  \sin(2\pi x_2) .
\] 
In this case, the exact solution of equation (\ref{2.4}) has the form
\[
 u(\bm x) = \nu_1^{-\alpha}\sin(\pi x_1)  \sin(\pi x_2) + \nu_2^{-\alpha} \sin(3\pi x_1)  \sin(2\pi x_2) ,
\]
\[
 \nu_1 = \pi^2,
 \quad \nu_2 = 13 \pi^2 .
\]
The error of the approximate solution was evaluated in norms $H$ and $H_A$:
\[
 \varepsilon = \|y - u\|,
 \quad \varepsilon_A = \|y - u\|_A .
\]
We use the uniform grid in space with $N_1 = N_2 = 100$ and the grid in time with $N = 20$ as the basic one
and the parameter $\alpha$ is equal to $0.5$, unless otherwise stated.

We start with numerical results for the problem (\ref{3.2}) obtained using (\ref{3.7}), (\ref{3.8}). 
The accuracy of the scheme with $\sigma = 1$ is presented in Table~\ref{tab-1}. Similar
data for the scheme (\ref{3.7}), (\ref{3.8} with $\sigma = 0.5$
are shown in Table~\ref{tab-2}. It is clear that the scheme of second-order accuracy ($\sigma = 0.5$)
provides good accuracy for a small number of time steps and demonstrates a weak dependence
on errors in specifying spectrum boundaries for the problem operator (compare data at the critical value 
$\theta = 1$ and $\theta = 0.5$). For the first-order scheme
($\sigma = 1$), the accuracy is significantly lower and influence of  $\theta$ is more strongly expressed.

\begin{table}
\begin{center}
 \caption{The error of the two-level scheme with weights for $\sigma = 1$}
 \begin{tabular}{cccccc}\label{tab-1}
  $N$  & 5 & 10 & 20 & 40 & 80  \\
  \hline
  $\varepsilon (\theta = 1)$     & 0.0120292 & 0.0066811 & 0.0035420 & 0.0018307 & 0.0009359 \\
  $\varepsilon_A (\theta = 1)$   & 0.1362146 & 0.0756551 & 0.0401086 & 0.0207297 & 0.0105974 \\
  $\varepsilon (\theta = 0.5)$   & 0.0255692 & 0.0148685 & 0.0081424 & 0.0042856 & 0.0022062 \\
  $\varepsilon_A (\theta = 0.5)$ & 0.2869894 & 0.1671826 & 0.0916372 & 0.0482531 & 0.0248449 \\
\end{tabular}
\end{center} 
\end{table} 

\begin{table}
\begin{center}
 \caption{The error of the two-level scheme with weights for $\sigma = 0.5$}
 \begin{tabular}{cccccc}\label{tab-2}
  $N$  & 5 & 10 & 20 & 40 & 80  \\
  \hline
  $\varepsilon (\theta = 1)$     & 0.0014351 & 0.0004130 & 0.0001236 & 0.0000484 & 0.0000296 \\
  $\varepsilon_A (\theta = 1)$   & 0.0162505 & 0.0046760 & 0.0013988 & 0.0005462 & 0.0003309 \\
  $\varepsilon (\theta = 0.5)$   & 0.0026070 & 0.0008303 & 0.0002392 & 0.0000720 & 0.0000287 \\
  $\varepsilon_A (\theta = 0.5)$ & 0.0295158 & 0.0094001 & 0.0027076 & 0.0008137 & 0.0003206 \\
\end{tabular}
\end{center} 
\end{table} 

The influence of a spacial grid size on the accuracy of a numerical solution 
for $\sigma = 1.0$ and $\sigma = 0.5$ is presented in Table~\ref{tab-3} for $\theta = 1$, $N = 80$.
In the case of $\sigma = 1.0$, there is no decreasing of the error with refining spatial grids, i.e.,
the main part of error results from the approximation in time. If $\sigma = 0.5$, then error
decreases with refining spatial grids. Table~\ref{tab-4} demonstrates numerical results obtained
with various values of the parameter $\alpha$ on the grid $N_1=N_2 = 100$, $N = 80$.

\begin{table}
\begin{center}
 \caption{Error for various grids in space}
 \begin{tabular}{cccccc}\label{tab-3}
  $N_1 = N_2$  & 25 & 50 & 100 & 200 & 400  \\
  \hline
  $\varepsilon (\sigma = 1)$     & 0.0010624 & 0.0009504 & 0.0009359 & 0.0009387 & 0.0009426 \\
  $\varepsilon_A (\sigma = 1)$   & 0.0119476 & 0.0107503 & 0.0105974 & 0.0106319 & 0.0106768 \\
  $\varepsilon (\sigma = 0.5)$   & 0.0002321 & 0.0000600 & 0.0000172 & 0.0000066 & 0.0000400 \\
  $\varepsilon_A (\sigma = 0.5)$ & 0.0025009 & 0.0006511 & 0.0001888 & 0.0000737 & 0.0000452 \\
\end{tabular}
\end{center} 
\end{table} 

\begin{table}
\begin{center}
 \caption{Error for various values of $\alpha$}
 \begin{tabular}{cccccc}\label{tab-4}
  $\alpha$  & 0.1 & 0.3 & 0.5 & 0.7 & 0.9  \\
  \hline
  $\varepsilon (\sigma = 1)$     & 0.0009628 & 0.0012876 & 0.0009359 & 0.0005621 & 0.0003062 \\
  $\varepsilon_A (\sigma = 1)$   & 0.0109020 & 0.0145807 & 0.0105974 & 0.0063651 & 0.0034674 \\
  $\varepsilon (\sigma = 0.5)$   & 0.0000925 & 0.0000277 & 0.0000172 & 0.0000090 & 0.0000045 \\
  $\varepsilon_A (\sigma = 0.5)$ & 0.0002772 & 0.0003097 & 0.0001888 & 0.0000950 & 0.0000433 \\
\end{tabular}
\end{center} 
\end{table} 

Possibilities of using splitting schemes for solving the problem (\ref{3.2}) 
are shown on the splitting with respect to spatial variables, where $p=2$ in (\ref{4.1})
and one-dimensional operators $A_1$ and $A_2$ are defined as it is given in Section \ref{sec:3}.
In our model problem with the Laplace operator, we have 
\[
 A_i \geq \delta_i I, 
 \quad \delta_i > 0, 
 \quad  i =1,2
\] 
with $\delta_i, \ i =1,2$ defined according to (\ref{5.1}).
For $D_i,  \ i =1,2$ (see (\ref{4.2})), we put
\[
 D_i = A_i - \chi_i I,
 \quad  \chi_i = \theta \delta_i ,
 \quad  i =1,2 .
\] 

In our case of the two-component splitting, the vector scheme (\ref{4.14}) takes the form
\[
\begin{split}
  \theta \delta  \frac{y_1^{n+1} - y_1^{n}}{\tau} & + 
  t^{n} D_1 \left ( \sigma_1  \frac{y_1^{n+1} - y_1^{n}}{\tau} 
  + (1-\sigma_1 ) \frac{y_1^{n} - y_1^{n-1}}{\tau} \right ) 
  + t^{n} D_2 \frac{y_2^{n} - y_2^{n-1}}{\tau } \\
  & + \alpha D_1 (\sigma_2 y_1^{n+1} + (1-\sigma_2) y_1^{n} )
   + \alpha D_2 y_2^{n} = 0,
\end{split}
\]
\[
\begin{split}
  \theta \delta  \frac{y_2^{n+1} - y_2^{n}}{\tau} & +
  t^{n} D_1 \frac{y_1^{n} - y_1^{n-1}}{\tau } +
  t^{n} D_2 \left ( \sigma_1  \frac{y_2^{n+1} - y_2^{n}}{\tau} 
  + (1-\sigma_1 ) \frac{y_2^{n} - y_2^{n-1}}{\tau} \right )  \\
  & + \alpha D_1 y_1^{n} + \alpha D_2 (\sigma_2 y_2^{n+1} + (1-\sigma_2) y_2^{n} )
   = 0 .
\end{split}
\]
For three-level schemes, we need to calculate separately  an approximate solution 
at the first level using a two-level scheme. 
Taking into account that we use the scheme with the first-order approximation in time,
the simplest explicit scheme can be applied. Then
\[
 \theta \delta  \frac{y_i^{1} - y_i^{0}}{\tau} + \alpha D y_i^{0} = 0,
 \quad i = 1,2 . 
\]
For our two-component splitting with $p=2$, we restrict ourselves (see Theorem~\ref{t-2}) 
to the case $\sigma_1 = \sigma_2 = 1$.

We can consider the error of individual components of the vector splitting scheme. 
Let
\[
 \varepsilon^{(i)} = \|y_i - u\|,
 \quad \varepsilon_A^{(i)} = \|y_i - u\|_A ,
 \quad i = 1,2 .
\]
In studying the splitting scheme, emphasis is on a dependence of error on a time step. 
We present the results of calculations by the splitting scheme
obtained on the basic grid in space $N_1 = N_2 = 100$ for $\alpha = 0.5$.
Table~\ref{tab-5} shows the accuracy of individual components of the vector splitting scheme
at the critical value of the parameter $\theta = 1$. Similar data for $\theta = 0.5$
are summarized in Table~\ref{tab-6}. Predictions demonstrate convergence
with refining grids in time. The accuracy of the splitting schemes is comparable with the accuracy 
of the implicit scheme without splitting (\ref{3.7}), (\ref{3.8}) at $\sigma = 1$. 

\begin{table}
\begin{center}
 \caption{The error of the splitting scheme with $\theta  = 1$}
 \begin{tabular}{cccccc}\label{tab-5}
  $N$  & 5 & 10 & 20 & 40 & 80  \\
  \hline
  $\varepsilon^{(1)}$     & 0.0084773 & 0.0032158 & 0.0012451 & 0.0005118 & 0.0002210 \\
  $\varepsilon_A^{(1)}$   & 0.0959950 & 0.0364151 & 0.0140991 & 0.0057948 & 0.0025025 \\
  $\varepsilon^{(2)}$     & 0.0251082 & 0.0080231 & 0.0029405 & 0.0012031 & 0.0005301 \\
  $\varepsilon_A^{(2)}$   & 0.2843178 & 0.0908516 & 0.0332978 & 0.0136229 & 0.0060020 \\
\end{tabular}
\end{center} 
\end{table} 

\begin{table}
\begin{center}
 \caption{The error of the splitting scheme with $\theta  = 0.5$}
 \begin{tabular}{cccccc}\label{tab-6}
  $N$  & 5 & 10 & 20 & 40 & 80  \\
  \hline
  $\varepsilon^{(1)}$     & 0.0063711 & 0.0045436 & 0.0019529 & 0.0007352 & 0.0002820 \\
  $\varepsilon_A^{(1)}$   & 0.0682378 & 0.0503844 & 0.0215774 & 0.0079984 & 0.0029949 \\
  $\varepsilon^{(2)}$     & 0.1398605 & 0.0399402 & 0.0130816 & 0.0049608 & 0.0020928 \\
  $\varepsilon_A^{(2)}$   & 1.5835471 & 0.4521733 & 0.1480800 & 0.0561468 & 0.0236847 \\
\end{tabular}
\end{center} 
\end{table}

\end{document}